\documentclass[12pt,psamsfonts,a4paper,reqno]{amsart}
\usepackage[left=1in,right=1in,bottom=1.2in,top=1.2in]{geometry}
\usepackage{tikz,amsthm,amsmath,amstext,amssymb,amscd,epsfig,euscript, mathrsfs, dsfont,pspicture,multicol, mathtools,graphpap,graphics,graphicx,times,enumerate,subfig,sidecap,wrapfig,color}
\usepackage{hyperref}
\hypersetup{
colorlinks,
citecolor=green,
filecolor=black,
linkcolor=blue,
urlcolor=blue,
}
\usepackage{bm}
\usepackage{amsmath}
\usepackage{mathrsfs}  
\usepackage{amsthm}
\usepackage{amssymb}
\usepackage{pgfplots}

\def \rr {{\mathbb R}}

\DeclareGraphicsExtensions{.png}
 \numberwithin{equation}{section}

\def\t{\tau}

\def\e{\eta}

\def\O{\Omega}

\def\S{\Sigma}

\newtheorem{definition}{Definition}[section]
\newtheorem{theorem}{Theorem}[section]
\newtheorem{proposition}{Proposition}[section]
\newtheorem{lemma}{Lemma}[section]

\newtheorem{remark}{Remark}[section]

	\title[Approximation of the $\delta$-shell interaction]{On the approximation of the $\delta$-shell interaction \\ for the 3-D Dirac operator.}
\author[Mahdi ZREIK]{Mahdi Zreik \textsuperscript{1}}
\address{\textsuperscript{1}Institut de Math\'ematiques de Bordeaux, UMR 5251, Universit\' e de Bordeaux 33405 Talence Cedex, FRANCE\\ and Departamento de Matem\'aticas, Universidad del Pa\' is Vasco, Barrio Sarriena s/n 48940 Leioa, SPAIN.}
	
\begin{document}
	\begin{abstract}
We consider the three-dimensional Dirac operator coupled with a combination of electrostatic and Lorentz scalar $\delta$-shell interactions. We approximate this operator with general local interactions $V$. Without any hypotheses of smallness on the potential $V$, converges in the strong resolvent sense to the Dirac Hamiltonian coupled with a $\delta$-shell potential supported on $\S$, a bounded smooth surface. However, the coupling constant depends nonlinearly on the potential $V.$
	\end{abstract}
\maketitle	
\tableofcontents
\section{Introduction}
Dirac Hamiltonians of the type $D_m + V$, where $V$ is a suitable perturbation, are used in many problems where the implications of special relativity play an important role. This is the case, for example, in the description of elementary particles such as quarks, or in the analysis of graphene, which is used in research for batteries, water filters, or photovoltaic cells. For these problems, mathematical investigations are still in their infancy. The present work studies the three-dimensional Dirac operator with a singular interaction on a closed surface $\Sigma $. \\\\
Mathematically, the Hamiltonian we are interested in can be formulated as follows
\begin{align}\label{DEL}
   \mathbb{D}_{\eta,\tau}=D_m+B_{\e,\t}\delta_{\S}= D_m + \big (\eta \,\mathbb{I}_{4} + \tau \beta\big)\delta_{\S},
\end{align}
where $B_{\e,\t}$ is the combination of the \textit{electrostatic} and \textit{Lorentz scalar} potentials of strengths $\eta$ and $\tau$, respectively. Physically, the Hamiltonian $\mathbb{D}_{\eta,\tau}$ is used as an idealized model for Dirac operators with strongly localized electric and massive potential near the interface $\S$ (\emph{e.g.}, an annulus), \emph{i.e.}, it replaces a Hamiltonian of the form 
\begin{align}\label{DEB}
\mathbb{H}_{\tilde{\eta},\tilde{\tau}}=D_m + B_{\tilde{\e},\tilde{\t}}= D_m  + \big (\tilde{\eta} \,\mathbb{I}_{4} + \tilde{\tau} \beta\big)\mathfrak{P}_{\S},
\end{align}
where $\mathfrak{P}_{\S}$ is a regular potential localized in a thin layer containing the interface $\S$. \\\\
The operators $\mathbb{D}_{\eta,\tau } $ have been studied in detail recently. Starting from the first directed work on spectral studies of Hamiltonian $\mathbb{D}_{\eta,\tau}$ in Ref. \cite{DES}, in which the authors treated the case that the surface is a sphere, assuming $\eta^{2} - \tau^{2}= -4.$ This is known as the \textit{confinement case} and in physics means the stability of a particle (\emph{e.g.}, an electron) in its initial region during time evolution,
\emph{i.e.}, if for time $t=0$ the particle is considered in a confined region $\Omega\subset\mathbb{R}^{3}$, then it cannot cross the surface $\partial\Omega$ to join the region $\mathbb{R}^{3}\setminus\overline{\Omega}$ for all $t > 0$.
Mathematically, this means that the Dirac operator under consideration decouples into a direct sum of two Dirac operators acting on $\Omega$ and $\mathbb{R}^{3}\setminus\overline{\Omega}$, respectively, with appropriate boundary conditions. After this work, spectral studies of Schrodinger's operators coupled with $\delta$-shell interaction flourished, while spectral studies of $\delta$-shell interaction of Dirac operators in deep stability were lifeless. \\\\ In 2014, the spectral studies of $\delta$-shell interaction of Dirac operators was revived in \cite{AMV1}, where the authors developed a new technique to characterize the self-adjointness of the free Dirac operator coupled to a $\delta$-shell potential. In a special case, they treated pure electrostatic $\delta$-shell interactions (\emph{i.e.}, $\tau = 0$) supported on the boundary of a bounded regular domain and proved that the perturbed operator is self-adjoint. The same authors continued the spectral study of the electrostatic case; for example, the existence of a point spectrum and related problems; see \cite{AMV2} and \cite{AMV3}.
\\\\
The approximation of the Dirac operator $\mathbb{D}_{\eta,\tau}$ by Dirac operators with regular potentials with shrinking support (\emph{i.e.}, of the form \eqref{DEB}) provides a justification of the considered idealized model. In the one-dimensional framework, the analysis is carried out in \cite{PS}, where Šeba showed that convergence in the sense of norm resolution is true. Subsequently, Hughes and Tušek prove strong resolvent convergence and norm resolvent convergence for Dirac operators with general point interactions in \cite{RJH1,RJH2} and \cite{Tus}, respectively. In 2D, \cite{CLMT} considered the approximation of Dirac operators with electrostatic, Lorentz scalar, and anomalous magnetic $\delta$-shell potentials on closed and bounded curves. Furthermore, in \cite{BHT} the authors examined the same question as in the paper \cite{CLMT}, but on a straight line. More precisely, taking parameters $(\tilde{\e},\tilde{\tau})\in\rr^2$ in \eqref{DEB} and a potential $\mathfrak{P}^{\varepsilon}_{\S}$ converging to $\delta_{\S}$ when $\varepsilon$ tends to $0$ (in the sense of distributions), then $D_m + \big (\tilde{\eta} \,\mathbb{I}_{4} + \tilde{\tau} \beta\big)\mathfrak{P}^{\varepsilon}_{\S}$ converges to the Dirac operator $\mathbb{D}_{\eta,\tau}$ with different coupling constants $(\e,\t)\in\rr^2$ which depend nonlinearly on the potential $\mathfrak{P}^{\varepsilon}_{\S}$.
\\\\
In the three-dimensional case, the situation seems to be even more complex, as recently shown in \cite{AF}. There, too, the authors were able to show the convergence in the norm resolvent sense in the non-confining case, however, a smallness assumption for the potential $\mathfrak{P}^{\varepsilon}_{\S}$ was required to achieve such a result. On the other hand, this assumption, unfortunately, prevents us from obtaining an approximation of the operator $\mathbb{D}_{\eta,\tau}$ with the parameters $\eta$ and $\tau$ which are more relevant from the physical or mathematical point of view. Believing this to be the case, the authors of the recent paper \cite{BHS} have studied and confirmed the approximation problem for two- and three-dimensional Dirac operators with delta-shell potential in norm resolvent sense. Without the smallness assumption of the potential $\mathfrak{P}^{\varepsilon}_{\S}$, no results could be obtained here either. Finally, we note that in the two- and three-dimensional setting a renormalization of the interaction strength was observed in \cite{CLMT,AF,BHS}.\\\\
Our main goal in this work is to develop new techniques that will allow us to establish the approximation in terms of the strong resolvent in the non-critical and non-confinement cases (\emph{i.e.}, when $\eta^2 - \tau^2 \neq \pm 4$) without the smallness assumed in \cite{AF} and to obtain results on how the initial parameters should be chosen so that the mathematical models reflect the physical reality in the correct way. \\ Let $m>0$, recall the free Dirac operator $D_m$ on $\rr^3$ defined by  $D_m:=- i \alpha \cdot\nabla + m\beta$, where   
\begin{align*}
	\alpha_k&=\begin{pmatrix}
	0 & \sigma_k\\
	\sigma_k & 0
	\end{pmatrix}\quad \text{for } k=1,2,3,
	\quad
	\beta=\begin{pmatrix}
	\mathbb{I}_{2} & 0\\
	0 & -\mathbb{I}_{2}
	\end{pmatrix},\quad 
\mathbb{I}_{2} :=\begin{pmatrix}
1& 0 \\
0 & 1
\end{pmatrix}, \\
&\text{and }\,
\sigma_1=\begin{pmatrix}
	0 & 1\\
	1 & 0
	\end{pmatrix},\quad \sigma_2=
	\begin{pmatrix}
	0 & -i\\
	i & 0
	\end{pmatrix} ,
	\quad
	\sigma_3=\begin{pmatrix}
	1 & 0\\
	0 & -1
	\end{pmatrix},
	\end{align*}  
 are the family of Dirac and Pauli matrices satisfy the anticommutation relations:
 \begin{align}\label{PMD}
   \lbrace \alpha_j,\alpha_k\rbrace = 2\delta_{jk}\mathbb{I}_4,\quad \lbrace \alpha_j,\beta\rbrace = 0,\quad \text{and} \quad \beta^{2} = \mathbb{I}_4,\quad j,k\in \lbrace 1,2,3\rbrace,  
 \end{align}
 where $\lbrace \cdot,\cdot\rbrace$ is the anticommutator bracket. We use the notation $\alpha \cdot x=\sum_{j=1}^{3}\alpha_j x_j$  for $x=(x_1,x_2,x_3)\in\rr^3$.  We recall that $(D_m, \mathrm{dom}(D_m))$ is self-adjoint (see, \emph{e.g.}, \cite[ Subsection~1.4]{Tha}), and that 
	\begin{align*}
	\mathrm{Sp}(D_m) =\mathrm{Sp}_{\mathrm{ess}}(D_m)=(-\infty,-m]\cup [m,+\infty).
	\end{align*}
Finally, we also give the Dirac operator coupled with  a combination of electrostatic, Lorentz scalar $\delta$-shell interactions of strength $\eta$ and $\tau$, respectively, which we will denote $D_{\eta,\tau}$ in what follows. Throughout this paper, for $\Omega\subset\mathbb{R}^3$ a bounded smooth domain with boundary $\Sigma:=\partial\Omega$, we refer to $H^{1}(\Omega,\mathbb{C}^{4}):=H^{1}(\O)^{4}$ as the first order Sobolev space 
\begin{align*}
\mathit{H}^1(\O)^4=\{ \varphi\in\mathit{L}^2(\O)^4: \text{ there exists } \tilde{\varphi}\in\mathit{H}^1(\rr^{3})^4 \text{ such that }  \tilde{\varphi}|_{\O} =\varphi\}.
\end{align*}
We denote by $H^{1/2}(\Sigma,\mathbb{C}^{4}):=H^{1/2}(\Sigma)^{4}$ the Sobolev space of order $1/2$ along the boundary $\Sigma$, and by $t_{\Sigma}: H^{1}(\Omega)^{4}\rightarrow H^{1/2}(\Sigma)^{4}$ the classical trace operator.
\begin{definition}\label{DiracelectroLorentz}  Let $\Omega$ be a bounded domain in $\mathbb{R}^3$ with a boundary $\S=\partial\Omega$. Let $(\eta, \tau)\in \mathbb{R}^{2}$. Then, $D_{\eta,\tau} =D_m + B_{\eta,\tau}\delta_{\Sigma} :=  D_m + (\eta \mathbb{I}_{4} + \tau\beta)\delta_{\Sigma} $ acting in $L^{2}(\rr^3)^{4}$ and defined as follows:  

\begin{align*}
&D_{\eta,\tau}f= D_m f_{+} \oplus D_m f_-, \text{ for all }f\in     \mathrm{dom}(D_{\eta,\tau}):=\lbrace f=f_{+} \oplus f_{-} \in H^{1}(\Omega)^{4} \oplus H^{1}(\mathbb{R}^{3}\setminus\overline{\Omega})^{4}: \\&\text{ the transmission condition (T.C) below holds in } H^{1/2}({\S})^{4}\rbrace.
\end{align*}
\begin{align}\label{DOMDL}
    \text{Transmission condition}: i\alpha\cdot\nu(t_{\S}f_{+} - t_{\S}f_-) + \dfrac{1}{2} (\eta \,\mathbb{I}_{4} + \tau\beta) (t_{\S}f_+ + t_{\S}f_-)=0,
\end{align}
\end{definition}
where $\nu$ is the outward pointing normal to $\Omega$.\qed\\\\
Recall that for $\eta^2 - \tau^2 \neq 0,4,$ the Dirac operator $(D_{\eta,\tau}, \mathrm{dom}(D_{\eta,\tau}))$ is self-adjoint and verifies the following assertions (see, e.g.,  \cite [Theorem 3.4, 4.1]{BEHL2})
\begin{itemize}
        \item [(i)] $\mathrm{Sp_{ess}}(D_{\eta,\tau})=(-\infty,m]\cup[m,+\infty).$
        \item [(ii)] $\mathrm{Sp_{dis}}(D_{\eta,\tau})$ is finite.\\
        \end{itemize}
 \textbf{Organization of the paper.} The present paper is structured as follows. We start with Section \ref{MMR}, where we define the model to be studied in our paper by introducing the family $\lbrace\mathscr{E}_{\eta,\tau,\varepsilon}\rbrace_{\varepsilon},$ which is the approximate Dirac operator family of operators $D_{\eta,\tau}$. We also discuss our main results by establishing Theorem \ref{maintheorem}. Moreover, in this section we give some geometric aspects characterizing the surface $\Sigma$, as well as some spectral properties of the Dirac operator coupled with the $\delta$-shell interaction presented in Lemma \ref{lemma3.1}. Section \ref{prooftheorem3.1} is devoted to the proof of Theorem \ref{maintheorem}, which approximates the Dirac operators with $\delta$-shell interaction by sequences of Dirac operators with regular potentials at the appropriate scale in the sense of strong resolvent. 
\section{Model and Main results}\label{MMR}
For a smooth bounded domain $\Omega\subset \mathbb{R}^{3}$, we consider an interaction supported on the boundary $\S:=\partial\Omega$ of $\Omega$.  The surface $\Sigma$ divides the Euclidean space into disjoint union $\mathbb{R}^{3}=\Omega_+ \cup\Sigma\cup\Omega_-$, where $\Omega_+:=\Omega$ is a bounded domain and $\Omega_-=\mathbb{R}^{3}\setminus\overline{\Omega_+}.$ We denote by $\nu$ and $\mathrm{dS}$ the unit outward pointing normal to $\Omega$ and the surface measure on $\Sigma$, respectively. 
We also denote by $f_{\pm}:=f\downharpoonright \Omega_{\pm}$ be the restriction of $f$ in $ \Omega_{\pm},$ for all $ C^{2}$–valued function $f$ defined on $\mathbb{R}^{3}.$ Then, we define the distribution $\delta_{\Sigma}f$ by
\begin{align*}
    \langle \delta_{\Sigma}f,\, g\rangle := \frac{1}{2} \int_{\S} (t_{\S} f_+ + t_{\S}f_-)\, g \, \mathrm{dS}, \quad \text{ for any test function}\quad g\in C^{\infty}_{0} (\mathbb{R}^{3})^{4},
    \end{align*}
    where $t_{\S}f_{\pm}$ is the classical trace operator defined below in Definition \ref{DiracelectroLorentz}.
Now, we explicitly construct regular symmetric potentials $V_{\eta,\tau,\varepsilon} \in L^{\infty}(\mathbb{R}^{3}; \mathbb{C}^{4\times4})$
supported on a tubular $\varepsilon$-neighbourhood of $\S$ and such that 
\begin{align*}
    V_{\eta,\tau,\varepsilon}  \xrightarrow[\varepsilon\to 0]{}  (\eta \,\mathbb{I}_{4} + \tau \beta)\delta_{\S} \quad \text{in the sense of distributions.}
\end{align*}
To explicitly describe the approximate potentials $V_{\eta,\tau,\varepsilon}$, we will introduce an additional notation. For $ \gamma>0,$ we define $\S_{\gamma}:= \lbrace x\in \mathbb{R}^{3}, \, \text{dist}(x,\S)< \gamma \rbrace$ a tubular neighborhood of $\S$ with width $\gamma$. For $\gamma>0$ small enough, $\S_{\gamma}$  is parametrized as 
\begin{align}\label{neighborhood}
    \S_{\gamma}= \lbrace x_{\S} + p\nu(x_{\S}),\, x_{\S}\in \S \quad \text{and} \quad p\in (-\gamma,\gamma)\rbrace.
\end{align}
For $0<\varepsilon<\gamma,$ let $h_{\varepsilon}(p):= \dfrac{1}{\varepsilon}h\left(\dfrac{p}{\varepsilon}\right)$, for all $p\in\mathbb{R}$,
with the function $h$ verifies the following 
\begin{align*}
     h\in L^{\infty}(\mathbb{R},\mathbb{R}),\quad
    \text{supp} \,h \subset (-1,1) \text{ and }\int_{-1}^{1}h(t)\,\text{d}t =1.
\end{align*}
Thus, we have:
\begin{align}\label{prophe}
    \text{supp} \,h_{\varepsilon} \subset (-\varepsilon,\varepsilon),\quad \int_{-\varepsilon}^{\varepsilon}h_{\varepsilon}(t)\,\text{d}t =1, \text{ and }\lim_{\varepsilon\rightarrow 0} h_{\varepsilon} = \delta_{0}\quad\text{in the sense of the distributions},
\end{align}
where $\delta_0$ is the Dirac $\delta$-function supported at the origin. Finally, for any $\varepsilon\in (0,\gamma)$, we define the symmetric approximate potentials $V_{\eta,\tau,\varepsilon}\in L^{\infty}(\mathbb{R}^{3},\mathbb{C}^{{4}\times{4}})$, as follows:
\begin{equation}\label{potential} 
	 V_{\eta,\tau,\varepsilon}(x):= 
	\left\{
	\begin{aligned}
 &B_{\eta,\tau}h_{\varepsilon}(p),\quad &\text{if}& \quad x=x_{\S} + p\nu(x_{\S})\in\S_{\gamma},\\
 &0,\quad  &\text{if}& \quad x\in \mathbb{R}^{3}\setminus \S_{\gamma}.
	\end{aligned}
	\right.
	\end{equation}
 It is easy to see that $\lim_{\varepsilon \rightarrow 0} V_{\eta,\tau,\varepsilon} = B_{\eta,\tau}\delta_{\S}$, in $\mathcal{D}^{'}(\mathbb{R}^{3})^{4}.$ For $0<\varepsilon<\gamma$, we define the family of Dirac operator $\lbrace\mathscr{E}_{\eta,\tau,\varepsilon}\rbrace_{\varepsilon}$ as follows:
 \begin{equation}\label{familleE}
\begin{aligned}
    &\mathrm{dom} (\mathscr{E}_{\eta,\tau,\varepsilon}):=\mathrm{dom} (D_m)= H^{1}(\mathbb{R}^{3})^{4},
   \\ & \mathscr{E}_{\eta,\tau,\varepsilon} \psi = D_{m} \psi + V_{\eta,\tau,\varepsilon}\psi,\quad\text{for all } \psi\in\mathrm{dom}(\mathscr{E}_{\eta,\tau,\varepsilon}).
\end{aligned}
 \end{equation}
The main purpose of the present paper is to study the strong resolvent limit of $\mathscr{E}_{\eta,\tau,\varepsilon}$ at $\varepsilon\rightarrow 0.$ To do this, we will introduce some notations and geometrical aspects which we will use in the rest of the paper.
\subsection{Notations and geometric aspects}\label{geo} 
Let $\S$ be parametrized by the family $\lbrace \phi_{j}, U_{j},V_{j},\rbrace_{j\in J}$ with $J$ a finite set, $U_j\subset \rr^{2},\,V_j \subset \rr^{3},\,\S\subset \bigcup_{j\in J}V_j$ and $\phi_j(U_{j})=V_j\cap\S\subset\S\subset\mathbb{R}^{3}$ for all $j\in J.$ We set $s=\phi_{j}^{-1}(x_{\S})$ for any $x_{\S}\in\S.$
\begin{definition}[Weingarten map]\label{Weingarten}
 For $x_{\S}=\phi_j (s)\in\S\cap V_j$ with $s\in U_j,$ one defines the Weingarten map (arising from the second fundamental form) as the following linear operator
\begin{align}
\begin{array}{rcl}
W_{x_{\S}}:=W(x_{\S}):T_{x_{\S}} &\to & T_{x_{\S}}\\
\partial_i \phi_j (s) &\mapsto & W(x_{\S})[\partial_i \phi_j] (s):=-\partial_i \nu (\phi_j (s)),
\end{array}
\end{align} where $T_{x_{\S}}$ denotes the tangent space of $\S$ on $x_{\S}$ and $\lbrace \partial_i \phi_j (s)\rbrace_{i=1,2}$ is a basis vector of $T_{x_{\S}}$.
\end{definition}
\begin{proposition}{\cite[Chapter 9 (Theorem 2), 12 (Theorem 2)]{JAT}.}\label{WP}
Let ${\S}$ be an $n-$surface in $\mathbb{R}^{n+1}$, oriented by the unit normal vector field $\nu$, and let $x\in {\S}$. Then, the Weingarten map verifies the following properties:
\begin{itemize}
    \item [(i)] Symmetric with respect to the inner product induced by the first fundamental form.
    \item [(ii)] Self-adjoint; that is $W_{x}(v)\cdot w=v \cdot W_{x}(w),$  for all $v,\,w\in T_{x}.$
    \item [(iii)] The eigenvalues $k_{1}(x),....,k_{n}(x)$ of the Weingarten map $W_{x}$ are called principal curvatures of $\S$ at $x$. Moreover, $k_{1}(x),....,k_{n}(x)$ uniformly bounded  on $\S$.
    \item [(iv)] The quadratic form associated with the Weingarten map at a point $x$ is called the second fundamental form of $\S$ at $x$.
\end{itemize}
\end{proposition}
The following theorem is the main result of this paper.
\begin{theorem}\label{maintheorem}
Let $(\eta,\tau)\in\mathbb{R}^{2}$, and denote by $d=\eta^2 - \tau^2$. Let $(\hat{\eta},\hat{\tau})\in\mathbb{R}^2$ be defined as follows:
\begin{equation}\label{parametres}
\begin{aligned}
&\bullet \text{if } d<0, \text{ then } 
    (\hat{\eta},\hat{\tau})=\dfrac{\mathrm{tanh}(\sqrt{-d}/2)}{(\sqrt{-d}/2)}(\eta,\tau),\\
       &\bullet \text{if } d=0, \text{ then } 
    (\hat{\eta},\hat{\tau})=(\eta,\tau),\\
       &\bullet \text{if } d>0 \text{ such that } d\neq (2k+1)^{2}\pi^2, \,k\in\mathbb{N}\cup\lbrace 0 \rbrace, \text{ then } 
    (\hat{\eta},\hat{\tau})=\dfrac{\mathrm{tan}(\sqrt{d}/2)}{(\sqrt{d}/2)}(\eta,\tau).\, 
\end{aligned}
\end{equation}
Now, let $\mathscr{E}_{\eta,\tau,\varepsilon}$ be defined as in \eqref{familleE} and $D_{\hat{\eta},\hat{\tau}}$ as in Definition \ref{DiracelectroLorentz}. Then, 
\begin{align}\label{relation}
    \mathscr{E}_{\eta,\tau,\varepsilon} \xrightarrow[\varepsilon\rightarrow 0]{} D_{\hat{\eta},\hat{\tau}} \quad \text{in the strong resolvent sense.}
\end{align}
\end{theorem}
\begin{remark}
    We mention that in this work we find approximations by regular potentials in the strong resolvent sense for the Dirac operator with $\delta$-shell potentials $\mathscr{E}_{\eta,\tau,\varepsilon}$ in the non-critical case (i.e., when $d \neq 4$) and non-confining case, (i.e., when $d\neq -4$) everywhere on $\Sigma.$ This is what we shall prove in the proof of Theorem \ref{maintheorem}.
\end{remark}
\subsubsection{Tubular neighborhood of $\S$}\label{tubular}
Recall that for $\Omega\subset\mathbb{R}^{3}$ a bounded domain with smooth boundary $\S$ parametrized by $\phi=\lbrace \phi_{j}\rbrace_{j\in J},$ we set $\nu_{\phi}=\nu\circ\phi: \S\longrightarrow \mathbb{R}^{3}$ the unit normal vector field which points outwards of $\O$, is independent of the particular choice of the positively oriented arc-length parametrization $\phi$.\\\\
For $\gamma>0,$ $\S_{\gamma}$ \eqref{neighborhood} is a tubular neighborhood of $\S$ of width $\gamma$.  We define the diffeomorphism $\Phi_{\phi} $ by:
\begin{align*}
    \Phi_{\phi}: U_{x_{\S}}\times (-\gamma,\gamma) &\longrightarrow \mathbb{R}^{3} \\ 
    (s,p) &\longmapsto \Phi_{\phi}(s,p) = \phi(s) + p\nu(\phi(s)). 
\end{align*}
For $\gamma$ be small enough, $\Phi_{\phi}$ is a smooth parametrization of $\S_{\gamma}$. Moreover, the matrix of the differential $ \mathrm{d}\Phi_{\phi}$ of $\Phi_{\phi}$ in the canonical basis of $\mathbb{R}^{3}$ is    
\begin{align}\label{gradfi}
    \mathrm{d}\Phi_{\phi}(s,p)= \Big( \partial_1 \phi(s) + p\,\mathrm{d}\nu(\partial_1 \phi)(s)\quad \partial_2 \phi(s) + p\,\mathrm{d}\nu(\partial_2 \phi)(s)\quad \nu_{\phi}(s)\Big).
\end{align}
Thus, the differential on $U_{x_{\Sigma}}$  and the differential on $(-\gamma,\gamma)$ of $\Phi_{\phi}$ are respectively given by 
\begin{equation}
\begin{aligned}
    &\mathrm{d}_{s}\Phi_{\phi}(s,p)= \partial_{i}\phi_j(s) - pW(x_{\S})\partial_{i}\phi_j(s) \quad \text{ for } i=1,2 \text{ and } x_{\S}\in\S,\\
    &\mathrm{d}_{p}\Phi_{\phi}(s,p)= \nu_{\phi}(s),
    \end{aligned}
    \end{equation}
    where $\partial_{i}\phi$, $\nu_{\phi}$ should be understood as column vectors, and $W(x_{\S})$ is the Weingarten map defined as in Definition \ref{Weingarten}.
    Next, we define 
    \begin{equation}\label{PP}
        \begin{aligned}
           & \mathscr{P}_{\phi}:= \Big(\Phi_{\phi}^{-1} \Big)_{1} : \S_{\gamma} \longrightarrow U_{x_{\S}}\subset\mathbb{R}^{2};\quad \mathscr{P}_{\phi}\big( \phi(s) + p\nu(\phi(s)) \big) = s\in\mathbb{R}^{2},\quad x_{\S}=\phi(s),\\
           & \mathscr{P}_{\perp}:= \Big(\Phi_{\phi}^{-1} \Big)_{2} : \S_{\gamma} \longrightarrow (-\gamma,\gamma);\quad \mathscr{P}_{\perp}\big( \phi(s) + p\nu(\phi(s)) \big) = p.
        \end{aligned}
    \end{equation}
    Using the inverse function theorem and thanks to \eqref{gradfi}, then we have for $x=\phi(s) + p\nu(\phi(s))\in\S_{\gamma}$ the following differential
    \begin{equation}\label{gradPP}
        \begin{aligned}
            \nabla \mathscr{P}_{\phi}(x)=\big ( 1-pW(s)\big)^{-1}t_{\phi}(s)\quad\text{and} \quad \nabla \mathscr{P}_{\perp}(x)=\nu_{\phi}(s),
        \end{aligned}
    \end{equation}
    with $t_{\phi}(s)=\partial_{i}\phi(s),$ $i=1,2.$
\subsection{Preparations for proof}
Before presenting the tools for the proof of Theorem \ref{maintheorem}, let us state some properties verified by the operator $D_{\eta,\tau}.$
 \begin{lemma}\label{lemma3.1}
     Let $(\eta,\tau)\in\mathbb{R}^{2}, $ and let $D_{\eta,\tau}$ be as in Definition \ref{DiracelectroLorentz}. Then, the following hold:
     \begin{itemize}
         \item [(i)] If $\eta^{2} - \tau^{2} \neq -4$, then there exists an invertible matrix $R_{\eta,\tau}$ such that a function $f = f_+ \oplus f_- \in H^1(\O_+)^4 \oplus H^1(\O_-)^4$ belongs to $\mathrm{dom}(D_{\eta,\tau})$ if and only if $t_{\S} f_{+} = R_{\eta,\tau} t_{\S}f_{-},$ with $R_{\eta,\tau}$ given by 
         \begin{align}\label{R}
            R_{\eta,\tau}:= \Big(\mathbb{I}_4  -  \dfrac{i\alpha\cdot\nu}{2}(\eta \,\mathbb{I}_{4} + \tau\beta)\Big)^{-1}\Big (\mathbb{I}_4 + \dfrac{i\alpha\cdot\nu}{2}(\eta \,\mathbb{I}_{4} + \tau\beta)\Big). 
         \end{align}
         \item[(ii)] If $\eta^{2} - \tau^{2} = -4,$ then  a function $f = f_+ \oplus f_- \in H^1(\O_+)^4 \oplus H^1(\O_-)^4$ belongs to $\mathrm{dom}(D_{\eta,\tau})$ if and only if 
\begin{align*}
\Big ( \,\mathbb{I}_{4} - \dfrac{i\alpha\cdot\nu}{2}(\eta \,\mathbb{I}_{4} + \beta\tau)\Big) t_{\S} f_{+} = 0 \quad \text{ and } \quad \Big ( \,\mathbb{I}_{4} + \dfrac{i\alpha\cdot\nu}{2}(\eta \,\mathbb{I}_{4} + \beta\tau)\Big) t_{\S} f_{-} = 0.
         \end{align*}
     \end{itemize}
 \end{lemma}
\textbf{Proof.} Using the transmission condition introduced in  \eqref{DOMDL}, then for assertion (i): for all $f=f_+ \oplus f_- \in \mathrm{dom}(D_{\eta,\tau}),$ we have that 
\begin{align*}
    \Big(i\alpha\cdot\nu +\dfrac{1}{2}(\eta \mathbb{I}_4 + \tau \beta)\Big)t_{\Sigma}f_+ =     \Big(i\alpha\cdot\nu  - \dfrac{1}{2}(\eta \mathbb{I}_4 + \tau \beta)\Big)t_{\Sigma}f_-.
\end{align*}
Thanks to properties in \eqref{PMD} and the fact that $(i\alpha\cdot\nu)^{-1} = -i\alpha\cdot\nu$, we get that
\begin{align}\label{relation}
    \Big( \,\mathbb{I}_{4} - M\Big) t_{\S} f_{+} = \Big( \,\mathbb{I}_{4} + M\Big) t_{\S} f_{-},
\end{align}
with $M$ a $4\times4$ matrix has the following form
\begin{align*}
    M=\dfrac{i\alpha\cdot\nu}{2}(\eta \,\mathbb{I}_{4} + \beta\tau),
\end{align*}
thus \eqref{R} is established.\\
Furthermore, as $d:=\eta^{2} - \tau^{2} \neq -4$ and $M^{2} = -\dfrac{d}{4} \mathbb{I}_4$,  $(\mathbb{I}_4 - M )(\mathbb{I}_4 + M) = \dfrac{4+d}{4}\mathbb{I}_4$, then $(\mathbb{I}_4 - M)$ is invertible and $(\mathbb{I}_4 - M)^{-1} = \dfrac{4}{4+d}(\mathbb{I}_4 + M)$. Consequently, using  \eqref{relation} we obtain that $t_{\Sigma}f_+ = R_{\eta,\tau}t_{\Sigma}f_-$ which $R_{\eta,\tau}$ has the following explicit form 
\begin{align}\label{ExplR}
    R_{\eta,\tau} = \dfrac{4}{4+d}\Bigg( \dfrac{4-d}{4}\mathbb{I}_4 + i\alpha\cdot\nu (\eta \mathbb{I}_4 + \tau\beta)\Bigg).
\end{align}
For assertion (ii), one just has to multiply \eqref{relation} by $(\mathbb{I}_4 \pm M)$ we get
\begin{align*}
    (\mathbb{I}_4 + M)^{2} t_{\Sigma}f_-= 0 \quad \text{and} \quad  (\mathbb{I}_4 - M)^{2} t_{\Sigma}f_+ = 0.
\end{align*}
This achieves the proof of Lemma \ref{lemma3.1}.\qed
\section{Proof of Theorem \ref{maintheorem}}\label{prooftheorem3.1}
Let $\lbrace\mathscr{E}_{\eta,\tau,\varepsilon}\rbrace_{\varepsilon \in (0,\gamma)}$ and $D_{\hat{\eta},\hat{\tau}}$ be as defined in \eqref{familleE} and Definition \ref{DiracelectroLorentz}, respectively. Since the singular interaction $V_{\eta,\tau,\varepsilon}$ are bounded and symmetric, then by the Kato-Rellich theorem, the operators $\mathscr{E}_{\eta,\tau,\varepsilon}$ are self-adjoint in $L^{2}(\mathbb{R}^{3})^{4}.$ Moreover, we know that $D_{\hat{\eta},\hat{\tau}}$ are self-adjoint and $\mathrm{dom}(D_{\hat{\eta},\hat{\tau}})\subset H^{1}(\mathbb{R}^{3}\setminus \Sigma)^{4}.$ Although the limiting operators and the limit operator are self-adjoint, it has been shown in \cite[Theorem VIII.26]{RS} that $\lbrace\mathscr{E}_{\eta,\tau,\varepsilon}\rbrace_{\varepsilon\in (0,\gamma)}$ converges in the strong resolvent sense to $D_{\hat{\eta},\hat{\tau}}$ as $\varepsilon\rightarrow 0$ if and only if it converges in the strong graph limit sense. The latter means that, for all $\psi \in \mathrm{dom}(D_{\hat{\eta},\hat{\tau}})$, there exists a family of vectors  $\lbrace \psi_{\varepsilon}\rbrace_{\varepsilon\in (0,\gamma)} \subset \mathrm{dom} (\mathscr{E}_{\eta,\tau,\varepsilon})$  such that 
\begin{align}\label{SGL}
   \mathrm{(a)} \,\lim_{\varepsilon\rightarrow 0} \psi_{\varepsilon} = \psi \quad\text{ and } \quad \mathrm{(b)}\,\lim_{\varepsilon\rightarrow 0} \mathscr{E}_{\eta,\tau,\varepsilon}\psi_{\varepsilon} = D_{\hat{\eta},\hat{\tau}}\psi \quad \text{in } L^{2}(\mathbb{R}^{3})^{4}.
\end{align}
Let $\psi\equiv \psi_+ \oplus \psi_- \in \mathrm{dom}(D_{\hat{\eta},\hat{\tau}}).$ From \eqref{parametres}, we have that 
\begin{align*}
    \hat{d} = \hat{\eta}^{2} - \hat{\tau}^{2} = -4\mathrm{tanh}^{2}(\sqrt{-d}/2),\quad\text{if}\quad d<0,\\
    \hat{d} = \hat{\eta}^{2} - \hat{\tau}^{2} = 4\mathrm{tan}^{2}(\sqrt{d}/2),\quad\text{if}\quad d>0,\\
    \hat{d} = \hat{\eta}^{2} - \hat{\tau}^{2} = 0,\quad\text{if}\quad d=0.
\end{align*}
In all cases, we have that $\hat{d} > -4$ (in particular $\hat{d}\neq -4)$. Then, by Lemma \ref{lemma3.1} (i), 
\begin{align}\label{}
    t_{\Sigma} \psi_+ = R_{\hat{\eta}, \hat{\tau}} t_{\Sigma} \psi_-,
\end{align}
where $R_{\hat{\eta}, \hat{\tau}}$ are given in \eqref{ExplR}.\\\\
Using the Definition \ref{DiracelectroLorentz}, we get that $t_{\Sigma}\psi_{\pm}\in H^{1/2}(\Sigma)^{4}.$\\\\
   $\bullet\,\,$ \textit{\textbf{Show that} }  \begin{align}\label{REXP}
       e^{i\alpha\cdot\nu B_{\eta,\tau}} = R_{\hat{\eta}, \hat{\tau}}.
    \end{align}
\\
Recall the definition of the family $\mathscr{E}_{\eta,\tau,\varepsilon}$ and $V_{\eta,\tau,\varepsilon}$ defined in \eqref{familleE} and \eqref{potential}, respectively. We have that 
\begin{align*}
    (i\alpha\cdot\nu B_{\eta,\tau})^{2}= (i\alpha\cdot\nu (\eta \mathbb{I}_4 + \tau\beta))^{2}=-(\eta^{2} - \tau^{2})=:D^{2},\quad \text{with } D= \sqrt{-(\eta^{2} - \tau^{2})}=\sqrt{-d}.
\end{align*}
Using this equality, we can write: $e^{i\alpha\cdot\nu B_{\eta,\tau}} = e^{-D} \Pi_- + e^{D}\Pi_+$, with $\pm D$ the eigenvalues of $ i\alpha\cdot\nu B_{\eta,\tau};$ and $\Pi_{\pm}$ the eigenprojections are given by:
\begin{align*}
    \Pi_{\pm}:= \dfrac{1}{2}\Bigg( \mathbb{I}_4 \pm \dfrac{i\alpha\cdot\nu B_{\eta,\tau}}{D}\Bigg).
\end{align*}
Therefore,
\begin{align*}
 e^{(i\alpha\cdot\nu B_{\eta,\tau})}&
 =\Bigg(\dfrac{ e^{D} + e^{-D}}{2} \Bigg)\mathbb{I}_4 + \dfrac{i\alpha\cdot\nu B_{\eta,\tau}}{D}\Bigg(\dfrac{e^{D} - e^{-D}}{2}\Bigg)\\&= \mathrm{cosh}(D) \mathbb{I}_4 + \dfrac{\mathrm{sinh}(D)}{D}(i\alpha\cdot\nu (\eta \mathbb{I}_4 + \tau\beta)).
\end{align*}
Now, the idea is to show \eqref{REXP}, \emph{i.e.}, it remains to show
\begin{align}\label{ERE}
    \dfrac{4}{4+\hat{d}}\Bigg( \dfrac{4-\hat{d}}{4} \mathbb{I}_4 + i\alpha\cdot\nu(\hat{\eta} \mathbb{I}_4 + \hat{\tau}\beta)\Bigg)-\mathrm{cosh}(D) \mathbb{I}_4 - \dfrac{\mathrm{sinh(D)}}{D}(i\alpha\cdot\nu (\eta \mathbb{I}_4 + \tau \beta))=0. 
\end{align}
To this end, set $\mathfrak{U}=\dfrac{4-\hat{d}}{4+\hat{d}}-\mathrm{cosh(D)}$ and $\mathfrak{V} = \dfrac{4}{4+\hat{d}} - \dfrac{\mathrm{sinh(D)}}{D}$. If we apply \eqref{ERE} to the unit vector $e_1 =(1\,\,0\,\,0\,\,0)^{t},$ then we get that $\mathfrak{U} = \mathfrak{V} = 0$. Hence, \eqref{ERE} makes sense  if and only if
\begin{align*}
 \quad \mathrm{cosh}(D) = \dfrac{4-\hat{d}}{4+\hat{d}} \quad\text{and} \quad \dfrac{\mathrm{sinh}(D)}{D}(\eta,\tau) = \dfrac{4}{4+\hat{d}}(\hat{\eta},\hat{\tau}).
\end{align*}
Consequently, we have $R_{\hat{\eta}, \hat{\tau}}= e^{i\alpha\cdot\nu B_{\eta,\tau}}$.\\\\
Moreover, dividing $\dfrac{\mathrm{sinh}(D)}{D}$ by $(1+\mathrm{cosh}(D))$ we get 
that     
\begin{align*}
     (\hat{\eta},\hat{\tau})= \dfrac{\mathrm{sinh}(D)}{1+\mathrm{cosh}(D)}\dfrac{1}{D/2}(\eta,\tau).
\end{align*}
Now, applying the elementary identity $\mathrm{tanh(\dfrac{\theta}{2})} = \dfrac{\mathrm{sinh}(\theta)}{1+\mathrm{cosh}(\theta)}$, for all $\theta\in\mathbb{C}\setminus\lbrace i(2k+1)\pi,\, k\in\mathbb{Z}\rbrace$.
We conclude that 
\begin{align*}
    \dfrac{\mathrm{tanh}(\sqrt{-d}/2)}{\sqrt{-d}/2}(\eta,\tau) = (\hat{\eta},\hat{\tau}),\quad \text{if  } d<0,
\end{align*}
and so, for $d>0$ we apply the elementary identity $-i\mathrm{tanh}(i\theta) = \mathrm{tan} (\theta)$ for all $ \theta\in \mathbb{C}\setminus \lbrace \pi(k+\dfrac{1}{2}),\, k\in \mathbb{Z}\rbrace$, then we get that
\begin{align*}
     \dfrac{\mathrm{tanh}(\sqrt{-d}/2)}{\sqrt{-d}/2}= \dfrac{\mathrm{tan}(\sqrt{d}/2)}{\sqrt{d}/2}.
\end{align*} 
Hence, we obtain that  $ \dfrac{\mathrm{tan}(\sqrt{d}/2)}{\sqrt{d}/2} (\eta,\tau) =(\hat{\eta},\hat{\tau}) \quad \text{if  } d>0$ such that $d\neq (2k+1)^{2}\pi^2$. Consequently, the equality $e^{i\alpha\cdot\nu B_{\eta,\tau}} = R_{\hat{\eta}, \hat{\tau}} $ is shown such that the following parameters verify:
\begin{equation}\label{base}
\begin{aligned}
    &\bullet \,\,\dfrac{\mathrm{tanh}(\sqrt{-d}/2)}{\sqrt{-d}/2}(\eta,\tau) = (\hat{\eta},\hat{\tau}),&\quad \text{if }& d<0,\\
   & \bullet\,\, \dfrac{\mathrm{tan}(\sqrt{d}/2)}{\sqrt{d}/2} (\eta,\tau) =(\hat{\eta},\hat{\tau}),& \quad \text{if }& d>0,\\
   & \bullet \,\, (\eta,\tau) = (\hat{\eta},\hat{\tau}),& \quad \text{if }& d=0.
\end{aligned}
\end{equation}
Moreover, the fact that $\int_{-\varepsilon}^{\varepsilon}h_{\varepsilon}(t)dt =1$ (see, \eqref{prophe}) with the statement \eqref{REXP} make it possible to write 
\begin{align}\label{1}
    \mathrm{exp}\Bigg[\big( -i\int_{-\varepsilon}^{0} h_{\varepsilon}(t) \,\mathrm{d}t\big)(\alpha\cdot\nu B_{\eta,\tau})\Bigg] t_{\Sigma}\psi_{+} =  \mathrm{exp}\Bigg[\big( i\int^{\varepsilon}_{0} h_{\varepsilon}(t) \,\mathrm{d}t\big)(\alpha\cdot\nu B_{\eta,\tau})\Bigg]t_{\Sigma}\psi_{-}.
\end{align}
$\bullet$ \textit{\textbf{Construction of the family $\boldsymbol{\lbrace\psi_{\varepsilon}\rbrace_{\varepsilon\in(0,\gamma)}}$.}}
For all $0<\varepsilon<\gamma$, we define the function $H_{\varepsilon}: \mathbb{R}\setminus\lbrace 0 \rbrace \rightarrow \mathbb{R}$ such that 
\begin{equation}\label{familleH} 
	 H_{\varepsilon}(p):= 
	\left\{
	\begin{aligned}
 &\int^{\varepsilon}_{p} h_{\varepsilon}(t) \,\mathrm{d}t,\quad &\text{if}& \quad 0<p<\varepsilon,\\
 & -\int_{-\varepsilon}^{p} h_{\varepsilon}(t) \,\mathrm{d}t,\quad  &\text{if}& \quad -\varepsilon<p<0,\\
 &0,\quad &\text{if}& \quad |p|\geq \varepsilon.
	\end{aligned}
	\right.
	\end{equation}
Clearly, $H_{\varepsilon}\in L^{\infty}(\mathbb{R})$ and supported in $(-\varepsilon,\varepsilon)$. The fact that $|| H_{\varepsilon}||_{L^{\infty}}\leq ||h||_{L^{1}}$, we get $\lbrace H_{\varepsilon}\rbrace_{\varepsilon}$ is
bounded uniformly in $\varepsilon.$ For all $\varepsilon\in(0,\gamma)$, the restrictions of $H_{\varepsilon}$ to $\mathbb{R}_{\pm}$ are uniformly continuous, so finite limits at $p=0$ exist, and differentiable a.e., with derivative being
bounded, since $h_{\varepsilon}\in L^{\infty}(\mathbb{R},\mathbb{R})$. Using these function, we set the matrix functions $\mathbb{U}_{\varepsilon}: \mathbb{R}^{3}\setminus \Sigma \rightarrow \mathbb{C}^{4\times4}$ such that
\begin{equation}\label{familleU} 
	 \mathbb{U}_{\varepsilon}(x):= 
	\left\{
	\begin{aligned}
 &e^{(i\alpha\cdot\nu) B_{\eta,\tau} H_{\varepsilon} (\mathscr{P}_{\perp}(x))},\quad &\text{if}& \quad x\in \Sigma_{\varepsilon}\setminus\Sigma,\\
 & \,\mathbb{I}_{4},\quad &\text{if}& \quad x\in\mathbb{R}^{3}\setminus\Sigma_{\varepsilon},
	\end{aligned}\quad \in L^{\infty}(\mathbb{R}^{3}, \mathbb{C}^{4\times4}),
	\right.
	\end{equation}
where the mappings $\mathscr{P}_{\perp}$ is defined as in \eqref{PP}. As the functions $\mathbb{U}_{\varepsilon}$ are bounded, uniformly in $\varepsilon$, and uniformly continuous in $\Omega_{\pm}$,
with a jump discontinuity across $\Sigma$, then $\forall x_{\Sigma}\in \Sigma$ and $y_{\pm}\in \Omega_{\pm}$, we get 
\begin{equation}\label{UEP}
\begin{aligned}
    \mathbb{U}_{\varepsilon}(x_{\Sigma}^{-})&:= \lim_{y_{-}\rightarrow x_{\Sigma}}\mathbb{U}_{\varepsilon}(y_-)= \mathrm{exp} \Big[i \Big(\int_{0}^{\varepsilon} h_{\varepsilon}(t) \,\mathrm{d}t \Big) (\alpha\cdot\nu (x_{\Sigma}))B_{\eta,\tau}\Big],\\
      \mathbb{U}_{\varepsilon}(x_{\Sigma}^{+})&:= \lim_{y_{+}\rightarrow x_{\Sigma}}\mathbb{U}_{\varepsilon}(y_+)= \mathrm{exp} \Big[ -i\Big(\int_{-\varepsilon}^{0} h_{\varepsilon}(t) \,\mathrm{d}t \Big) (\alpha\cdot\nu (x_{\Sigma}))B_{\eta,\tau}\Big].
\end{aligned}
\end{equation}
Thus, we construct $\psi_{\varepsilon}$ by $\psi_{\varepsilon} = \psi_{\varepsilon,+} \oplus \psi_{\varepsilon,-} := \mathbb{U}_{\varepsilon} \psi\in L^{2} (\mathbb{R}^{3})^{4}.$\\\\
Since $\mathbb{U}_{\varepsilon}$ are bounded, uniformly in $\varepsilon$, using the construction of $\psi_{\varepsilon}$ we get that $\psi_{\varepsilon} - \psi:= (\mathbb{U}_{\varepsilon} -\mathbb{I}_4 )\psi$. Then, by the dominated convergence theorem and the fact that $\mathrm{supp}\, (\mathbb{U}_{\varepsilon} -\mathbb{I}_4 )\subset | \Sigma_{\varepsilon}|$ with $| \Sigma_{\varepsilon}|\rightarrow 0$ as $\varepsilon \rightarrow 0$, it is easy to show that 
\begin{align}\label{a}
    \psi_{\varepsilon}\xrightarrow[\varepsilon\to 0]{}\psi\quad\text{in } L^{2}(\mathbb{R}^{3})^{4}. 
\end{align}
This achieves assertion (a). \\\\
$\bullet$ \textit{\textbf{Show that} $\boldsymbol{\psi_{\varepsilon} \in \mathrm{dom}(\mathscr{E}_{\eta,\tau,\varepsilon})= H^{1}(\mathbb{R}^{3})^{4}}$.}
This means that we must show, for all $0<\varepsilon<\gamma$, 
\begin{align*}
   \mathrm{(i)} \,\psi_{\varepsilon,\pm}\in H^{1}(\Omega_{\pm})^{4} \quad \text{and} \quad  \mathrm{(ii)}\,t_{\Sigma}\psi_{\varepsilon,+} = t_{\Sigma}\psi_{\varepsilon,-} \in H^{1/2}(\Sigma)^{4}.
\end{align*}
Let us show point (i). By construction of $\psi_{\varepsilon}$, we have $\psi_{\varepsilon} \in L^{2}(\mathbb{R}^{3})^{4}.$ It remains  to have $\partial_{j}\mathbb{U}_{\varepsilon} \in L^{2}(\mathbb{R}^{3})^{4}$, for $j=1,2,3.$ To do so, recall the parametrization $\phi: U \rightarrow \Sigma\subset\mathbb{R}^{3}$ of $\Sigma$ defined at the beginning of part \eqref{geo} and let $A$ a $4\times 4$ matrix such that $A(s):=i\alpha\cdot\nu(\phi(s)) B_{\eta,\tau},$ for $s=(s_1,s_2)\in U\subset\mathbb{R}^{2}.$ Thus, the matrix functions $\mathbb{U}_{\varepsilon}$ in \eqref{familleU} can be written 
\begin{equation}\label{familleU1} 
	 \mathbb{U}_{\varepsilon}(x)= 
	\left\{
	\begin{aligned}
 &e^{A(\mathscr{P}_{\phi}(x)) H_{\varepsilon} (\mathscr{P}_{\perp}(x))},\quad &\text{if}& \quad x\in \Sigma_{\varepsilon}\setminus\Sigma,\\
 & \,\mathbb{I}_{4},\quad &\text{if}& \quad x\in\mathbb{R}^{3}\setminus\Sigma_{\varepsilon},
	\end{aligned}\quad \in L^{\infty}(\mathbb{R}^{3}, \mathbb{C}^{4\times4}),
	\right.
	\end{equation}
where $\mathscr{P}_{\phi}$ is defined as in \eqref{PP}.\\\\
For $j=1,2,3,$ $\mathrm{supp}\,\partial_{j}\mathbb{U}_{\varepsilon}\subset \Sigma_{\varepsilon}$. Furthermore, it was mentioned in \cite[Eq.(4.1)]{RMW} that for all $x\in \Sigma_{\varepsilon}\setminus\Sigma$, $\partial_{j}\mathbb{U}_{\varepsilon}$ can be written as follows
\begin{equation}\label{djU}
\begin{aligned}
\partial_{j}\mathbb{U}_{\varepsilon}(x)=\int_{0}^{1} \Bigg[\mathrm{exp} \Big( z A(\mathscr{P}_{\phi}(x)) H_{\varepsilon} (\mathscr{P}_{\perp}(x)) \Big)\partial_{j}\Big(A(\mathscr{P}_{\phi}(x)) H_{\varepsilon} (\mathscr{P}_{\perp}(x))\Big)\times \\ \mathrm{exp} \Big( (1-z) A(\mathscr{P}_{\phi}(x)) H_{\varepsilon} (\mathscr{P}_{\perp}(x)) \Big)\Bigg]\mathrm{d}z.
\end{aligned}
\end{equation}
Let $x=\phi(s) + p \nu(\phi(s))\in \Sigma_{\gamma},$ and recall the definition of the mappings $\mathscr{P}_{\phi}(x)$ and $\mathscr{P}_{\perp}(x)$ introduced in \eqref{PP}. Based on the quantities \eqref{gradPP} (with $s=\mathscr{P}_{\phi}(x)$ and $p=\mathscr{P}_{\perp}(x)$), we get that
\begin{align}\label{djU2}
    \partial_{j}\Big(A(\mathscr{P}_{\phi}(x)) H_{\varepsilon} (\mathscr{P}_{\perp}(x)) \Big) = \partial_{s}A(s)(1-p W(s))^{-1} (t_{\phi}(s))_{j}H_{\varepsilon}(p) - A(s) h_{\varepsilon}(p)(\nu_{\phi}(s))_{j}.
\end{align} 
Therefore, $\partial_{j}\mathbb{U}_{\varepsilon}$ has the following form 
\begin{equation}\label{djU3}  
\begin{aligned}
\partial_{j}\mathbb{U}_{\varepsilon}(x)&= - A(s) h_{\varepsilon}(p)(\nu_{\phi}(s))_{j}\mathbb{U}_{\varepsilon}(x)\\& + \int_{0}^{1} e^{z A(s) H_{\varepsilon} (p)}\big[\partial_{s}A(s)(1-p W(s))^{-1} (t_{\phi}(s))_{j}H_{\varepsilon}(p)\big]e^{(1-z) A(s) H_{\varepsilon} (p)}\,\mathrm{d}z.
\end{aligned}
\end{equation}
Set by $\mathbb{E}_{\varepsilon,j}$ the second term of the right part of equality \eqref{djU3}, \emph{i.e.},
\begin{align}\label{EE}
\mathbb{E}_{\varepsilon,j} = \int_{0}^{1} e^{z A(s) H_{\varepsilon} (p)}\big[\partial_{s}A(s)(1-p W(s))^{-1} (t_{\phi}(s))_{j}H_{\varepsilon}(p)\big]e^{(1-z) A(s) H_{\varepsilon} (p)}\,\mathrm{d}z.
\end{align}
Then, thanks to the third property of the Proposition \ref{WP} verified by the Weingarten map, the matrix-valued functions $\mathbb{E}_{\varepsilon,j}$ are bounded, uniformly for $0<\varepsilon<\gamma$, and $\mathrm{supp}\,\mathbb{E}_{\varepsilon,j}\subset \Sigma_{\varepsilon}$. Moreover, we have $\mathbb{U}_\varepsilon$ and $\partial_j\mathbb{U}_\varepsilon \in L^{\infty}(\Omega_{\pm}, \mathbb{C}^{4\times4})$. Hence, for all $\psi_{\pm}\in H^{1}(\Omega_{\pm})^{4}$ we have that $\psi_{\varepsilon, \pm}= \mathbb{U}_{\varepsilon}\psi_\pm\in H^{1}(\Omega_{\pm})^{4}$ and statement (i) is proved.   \\\\
Now, we show point (ii). As $\psi_{\varepsilon, \pm} \in H^{1}(\Omega_\pm)^{4}$,  we get that $ t_{\Sigma}\psi_{\varepsilon,\pm}\in H^{1/2}(\Sigma)^{4}.$ On the other hand, it have been showed in \cite[Chapter 4 (p.133)]{EG}, for a.e., $x_{\Sigma}\in\Sigma$ and $r>0$, that 
\begin{align*}
    t_{\Sigma}\psi_{\varepsilon,\pm} (x_{\Sigma}) = \lim_{r\rightarrow 0 }\dfrac{1}{| B(x_{\Sigma},r)|}\int_{\Omega_{\pm}\cap B(x_{\Sigma},r)}\psi_{\varepsilon}(y)\,\mathrm{d}y= \lim_{r\rightarrow 0 }\dfrac{1}{| B(x_{\Sigma},r)|}\int_{\Omega_{\pm}\cap B(x_{\Sigma},r)}\mathbb{U}_{\varepsilon}(y)\psi(y)\,\mathrm{d}y,
\end{align*}
and so, similarly, 
\begin{align*}
    \mathbb{U}_{\varepsilon}(x_{\Sigma}^{\pm})t_{\Sigma}\psi_{\pm} (x_{\Sigma}) = \lim_{r\rightarrow 0 }\dfrac{1}{| B(x_{\Sigma},r)|}\int_{\Omega_{\pm}\cap B(x_{\Sigma},r)}\mathbb{U}_{\varepsilon}(x_{\Sigma}^{\pm})\psi(y )\,\mathrm{d}y.
\end{align*}
As $\mathbb{U}_{\varepsilon}$ is continuous in $\overline{\Omega_{\pm}}$, we get $ t_{\Sigma}\psi_{\varepsilon,\pm} (x_{\Sigma})=\mathbb{U}_{\varepsilon}(x_{\Sigma}^{\pm})t_{\Sigma}\psi_{\pm} (x_{\Sigma}).$
Consequently, \eqref{1} with  \eqref{UEP} give us that $t_{\Sigma} \psi_{\varepsilon,+} = t_{\Sigma} \psi_{\varepsilon,-} \in H^{1/2}(\Sigma)^{4}$. With this, (ii) is valid and $\psi_{\varepsilon}\in \mathrm{dom}(\mathscr{E}_{\eta,\tau,\varepsilon})$.\\\\
To complete the proof of Theorem \ref{maintheorem}, it remains to show the property (b), mentioned in \eqref{SGL}. Since $ (\mathscr{E}_{\eta,\tau,\varepsilon}\psi_{\varepsilon} - D_{\hat{\eta}, \hat{\tau}}\psi)$ belongs to $L^{2}(\mathbb{R}^{3})^{4}$, it suffices to prove the following:
\begin{align}\label{EDO}
     \mathscr{E}_{\eta,\tau,\varepsilon}\psi_{\varepsilon,\pm} - D_{\hat{\eta}, \hat{\tau}}\psi_{\pm} \xrightarrow[\varepsilon\to 0]{}0\quad \text{in }  L^{2}(\Omega_{\pm})^{4}.
\end{align}
To do this, let $\psi\equiv\psi_+\oplus\psi_- \in \mathrm{dom} (D_{\hat{\eta}, \hat{\tau}})$ and $\psi_{\varepsilon}\equiv\psi_{\varepsilon,+}\oplus\psi_{\varepsilon,-}\in\mathrm{dom} (\mathscr{E}_{\eta,\tau,\varepsilon})$. We have
\begin{equation}\label{diff}
\begin{aligned}
    \mathscr{E}_{\eta,\tau,\varepsilon}\psi_{\varepsilon,\pm} - D_{\hat{\eta}, \hat{\tau}}\psi_{\pm} &= -i\alpha\cdot\nabla\psi_{\varepsilon,\pm}+ m\beta\psi_{\varepsilon,\pm} + V_{\eta,\tau,\varepsilon}\psi_{\varepsilon,\pm} + i\alpha\cdot\nabla\psi_\pm - m\beta\psi_\pm \\&
    = -i\alpha\cdot\nabla(\mathbb{U}_{\varepsilon}\psi_\pm) + i\alpha\cdot\nabla\psi_\pm + m\beta (\mathbb{U}_{\varepsilon} - \mathbb{I}_4)\psi_\pm + V_{\eta,\tau,\varepsilon}\psi_{\varepsilon,\pm}\\& = -i\sum_{j=1}^{3}\alpha_j\big[(\partial_{j}\mathbb{U}_{\varepsilon})\psi_\pm + (\mathbb{U}_{\varepsilon} - \mathbb{I}_4)\partial_{j}\psi_\pm\big] + m\beta (\mathbb{U}_{\varepsilon} - \mathbb{I}_4)\psi_\pm + V_{\eta,\tau,\varepsilon}\psi_{\varepsilon,\pm}.
\end{aligned}
\end{equation}
Using the form of $\partial_{j}\mathbb{U}_{\varepsilon}$ given in \eqref{djU3},  the quantity $-i\sum_{j=1}^{3}\alpha_j(\partial_{j}\mathbb{U}_{\varepsilon})\psi_\pm$ yields 
\begin{align*}
    -i\sum_{j=1}^{3}\alpha_j(\partial_{j}\mathbb{U}_{\varepsilon})\psi_\pm &= -i\sum_{j=1}^{3} \alpha_{j} \big[ -i\alpha\cdot\nu V_{\eta,\tau,\varepsilon}\nu_{j} \mathbb{U}_{\varepsilon}\psi_\pm + \mathbb{E}_{\varepsilon,j}\psi_\pm\big] \\& = -(\alpha\cdot\nu)^{2} V_{\eta,\tau,\varepsilon}\psi_{\varepsilon,\pm} - i\sum_{j=1}^{3}\alpha_j\mathbb{E}_{\varepsilon,j}\psi_\pm\\& = -V_{\eta,\tau,\varepsilon}\psi_{\varepsilon,\pm} + \mathbb{R}_{\varepsilon}\psi_\pm,
\end{align*}
where $ \mathbb{E}_{\varepsilon,j}$ is given in \eqref{EE} and  $\mathbb{R}_{\varepsilon}= - i\sum_{j=1}^{3}\alpha_j\mathbb{E}_{\varepsilon,j},$ a matrix-valued functions in $L^{\infty}(\mathbb{R}^{3}, \mathbb{C}^{4\times4})$, verifies the same property of $\mathbb{E}_{\varepsilon,j}$  given in \eqref{EE}, for $\varepsilon\in(0,\gamma).$ Thus, \eqref{diff} becomes 
\begin{align*}
    \mathscr{E}_{\eta,\tau,\varepsilon}\psi_{\varepsilon,\pm} - D_{\hat{\eta}, \hat{\tau}}\psi_\pm= -i\sum_{j=1}^{3}\alpha_j\big[(\mathbb{U}_{\varepsilon} - \mathbb{I}_4)\partial_{j}\psi_\pm\big] + m\beta (\mathbb{U}_{\varepsilon} - \mathbb{I}_4)\psi_\pm + \mathbb{R}_{\varepsilon}\psi.
\end{align*}
 Since $\psi_{\pm}\in H^{1}(\Omega_{\pm})^{4}$, $(\mathbb{U}_{\varepsilon} - \mathbb{I}_4)$ and $ \mathbb{R}_{\varepsilon}$ are bounded,  uniformly in  $\varepsilon\in(0,\gamma)$ and supported in $\Sigma_{\varepsilon},$ and $|\Sigma_{\varepsilon}|$ tends to 0 as $\varepsilon\rightarrow 0$. By the dominated convergence theorem, we conclude that 
\begin{align}\label{relfin}
    \mathscr{E}_{\eta,\tau,\varepsilon}\psi_{\varepsilon,\pm} - D_{\hat{\eta}, \hat{\tau}}\psi_\pm \xrightarrow[\varepsilon\rightarrow 0 ]{} 0,\quad \text{holds in } L^{2}(\Omega_\pm)^{4},
\end{align}
and this achieves the assertion \eqref{EDO}. \\\\
Thus, the two conditions mentioned in \eqref{SGL} (\emph{i.e.}, (a) and (b)) of the convergence in the strong graph limit sense are proved (see, \eqref{a} and \eqref{relfin}). Also, note that the latter remains stable with respect to bounded symmetric perturbations (in our case $m\beta(\mathbb{U}_{\varepsilon} - \mathbb{I}_4 )$ with $m>0$, so we can assume $m=0$). Hence, the family $\lbrace\mathscr{E}_{\varepsilon}\rbrace_{\varepsilon\in (0,\gamma)}$ converges in the strong resolvent sense to $D_{\hat{\eta}, \hat{\tau}}$ as $\varepsilon\rightarrow 0.$ The proof of the Theorem \ref{maintheorem} is complete.
\qed
\section{Acknowledgment}
I wish to express my gratitude to my thesis advisor Luis VEGA for suggesting the problem and for many stimulating conversations, for his patient advice, and enthusiastic encouragement. I would also like to thank my supervisor Vincent BRUNEAU for his reading, comments, and helpful criticism.


\end{document}